\documentclass[a4paper,10pt]{article}
\usepackage{amssymb,amsmath}
\newcommand{\Q}{\mathbb{Q}}
\newcommand{\R}{\mathbb{R}}
\newcommand{\C}{\mathbb{C}}
\newcommand{\N}{\mathbb{N}}
\newcommand{\Proj}{\mathbb{P}}

\newcommand{\x}{{\bf x}}

\newcommand{\Z}{\mathbb{Z}}
\newcommand{\Qb}{\overline{\Q}}
\newcommand{\cN}{\mathcal{N}}

\newcommand{\ep}{\varepsilon}
\renewcommand{\mod}[1]{\hspace{-2.9mm}\pmod{#1}}

\newtheorem{theorem}{Theorem}
\newtheorem{lemma}{Lemma}
\newtheorem{prop}{Proposition}

\begin{document}
\title{Sums and Differences of Three $k$-th Powers}
\author{D.R. Heath-Brown\\Mathematical Institute, Oxford}
\date{}
\maketitle

\section{Introduction}

If $k\ge 2$ is a positive integer the number of representations of a
positive integer $N$ as either $x_1^k+x_2^k=N$ or $x_1^k-x_2^k=N$,
with integers $x_1$ and $x_2$, is
finite.  Moreover it is easily shown to be $O_{\ep}(N^{\ep})$, for any 
$\ep>0$.  It is known that
if $k=2$ or 3 then the number of representations is unbounded as $N$
varies, but it is conjectured that the number of representations is
bounded for $k\ge 4$.  Indeed for $k\ge 5$ we know of no $N$ for which
there are two or more essentially different representations.

This paper is primarily concerned with the analogous questions when one
has three $k$-th powers.  When $k=2$ or $3$ the equation
$x_1^k+x_2^k-x_3^k=N$ may have infinitely many solutions, as the
identities
\[(2t-1)^2+(t^2-t-1)^2-(t^2-t+1)^2=1\]
and
\[(6t^3+1)^3+(-6t^3+1)^3-(6t^2)^3=2\]
show.  Thus it is natural to count solutions with restrictions on the size of
the $x_i$.
If we require $\max |x_i|\le B$, say, then there will trivially be
$O_{\ep}(B^{1+\ep})$ solutions, since there are $O_{\ep}(B^{\ep})$
solutions for each value of $x_1^k\not=N$.  However, since the equation
$x_1^k+x_2^k-x_3^k=1$ has trivial solutions $(x_1,x_,x_3)=(1,t,t)$ it
is not possible to improve this bound beyond $O(B)$ in general.  We
shall say that the equation $F(x_1,x_2,x_3)=N$ has a parametric solution 
of degree $d$, if there are
polynomials $f_1(t),f_2(t),f_3(t)\in\Z[t]$ of maximum degree $d$, such
that $F(f_1(t),f_2(t),f_3(t))=N$, identically in $t$.  We write $S_d$
for the set of solutions given by parameterizations of degree at most
$d$, and we define
\[\cN(B)=\cN(B;N,F,d)=\#\{\x\in\Z^3:F(\x)=N,\, \max|x_i|\le B,\,
\x\not\in S_d\}.\]
We will consider arbitrary integral forms $F(x_1,x_2,x_3)$ of degree
$k$, although our main interest will be in diagonal forms $x_1^k\pm
x_2^k\pm x_3^k$. For these latter forms we shall refer to the trivial
solutions in which one of the individual terms $x_1^k, \pm x_2^k$ or $ \pm
x_3^k$ is equal to $N$, as being ``special solutions''.

The problem above has been investigated by the author \cite[Theorem
13]{annals} and \cite[Theorem 17]{cime}, where it is shown that
\begin{equation}\label{first}
\cN(B)\ll_{\ep,k}B^{2/\sqrt{k}+2/(k-1)+\ep}
\end{equation}
for any $\ep>0$, in the case $F(\x)=x_1^k+x_2^k+x_3^k$.  The exponent 
is non-trivial, that is to say we have $2/\sqrt{k}+2/(k-1)<1$, 
for $k\ge 8$.  However work by Salberger \cite{salb},
combined with methods from the papers above, allows one to replace the
exponent by $2/\sqrt{k}+\ep$, which is non-trivial for $k>4$.  Indeed
the result extends to the forms $F(\x)=x_1^k\pm x_2^k\pm x_3^k$.

It is a nice feature of these results that they hold completely
uniformly in
$N$.  However the goal of the present paper is to show how one can
obtain considerably sharper exponents in our estimate for $\cN(B)$
providing that one accepts a bound which is not completely 
uniform in $N$.  More
precisely we will prove the following theorems.
\begin{theorem}\label{th9/10}
Let $F(x_1,x_2,x_3)\in\Z[x_1,x_2,x_3]$ be a non-singular form of 
degree at least $3$, and let $\ep>0$ be given.  
Then if $N\ll_F B^{3/13}$ is a natural number, 
all integer points for which
\[|F(x_1,x_2,x_3)|\le N,\;\;\;\; B/2<\max|x_i|\le B\]
lie on a union of $O_F(B^{9/10}N^{1/10})$ plane projective conics
$C_i(x_1,x_2,x_3)=0$, with $C_i(x_1,x_2,x_3)\in\Z[x_1,x_2,x_3]$.

For such $N$ we have $\cN(B;N,F,1)=O_{F,\ep}(B^{9/10+\ep}N^{1/10})$.  
The number
of essentially different linear parameterizations is
bounded uniformly in terms of the degree of $F$.
Moreover 
in the cases when $F(x_1,x_2,x_3)=x_1^k\pm x_2^k\pm x_3^k$ there are 
$O_{\ep}(B^{9/10+\ep}N^{1/10})$ solutions apart from any special solutions.
\end{theorem}

\begin{theorem}\label{th1/k}
Let $F(x_1,x_2,x_3)\in\Z[x_1,x_2,x_3]$ be a non-singular form of 
degree $k\ge 3$, and let $N\ll_F B$ be a natural number.  Then 
\[\cN(B;N,F,[\frac{k}{10}])\ll_{F}B^{10/k}.  \]
The number
of essentially different parameterizations of degree at most $[k/10]$ is
bounded in terms of $k$ alone.  Moreover 
in the cases when 
\[F(x_1,x_2,x_3)=x_1^k\pm x_2^k\pm x_3^k\]
there are 
$O_{k}(B^{10/k})$ solutions apart from any special solutions.
\end{theorem}

Thus we get an the exponent $9/10$ for $B$, which is  non-trivial 
for every $k\ge 3$.
Moreover we get a general exponent for $B$ of order $1/k$, where previously
the best exponent had order $1\sqrt{k}$.
It is clear that the exponent in Theorem \ref{th1/k} might
be improved slightly.  However we have tried to give an argument which
leads to an exponent of order $1/k$ in as simple a manner as possible.
It would also be possible to establish a version of Theorem
\ref{th1/k} in which the admissible range for $N$ was extended.  In particular
one could still achieve an exponent for $B$ which was of order $1/k$,
while allowing any $N\le B^{k/2}$, for example. In fact the results
above are just two examples of a range of possible estimates with
different exponents for $B$ and $N$, depending on the degree of $F$ and
the size of $N$ relative to $B$.

We also observe that rather little use is made of our assumption that
$F$ should be non-singular.  It seems likely that related results
could be obtained without this condition, and indeed our treatment of
Theorem \ref{sqfr} below corresponds to the use of $F(x,y,z)=y^{k-1}z-x^k$.

Our results should be seen as examples of the ``determinant method''
developed by the author in \cite{annals}.  The idea has its origins in
the work of Bombieri and Pila \cite{BP}, and also of Elkies
\cite{Elk}.  Bombieri and Pila used a real-variable method to
investigate integral points on plane curves.  Their curves were
analytic, but not necessarily algebraic. Elkies, who also used a real
variable method, was more interested in
search algorithms for integral points on or near an algebraic
curve. The author replaced the real-variable approach with a $p$-adic
method, which simplifies many of the arguments, and extended the scope
to include rational and integral points on affine and projective
algebraic varieties of arbitrary dimension.  Thus, for example, the
bound (\ref{first}) arises from consideration of integral points on
the affine surface $x_1^k+x_2^k+x_3^k=N$.  In contrast, our current
work, in which we think of $N$ as being relatively small compared to
$B$, approaches the problem by examining points near to the
projective curve $x_1^k+x_2^k+x_3^k=0$.  Loosely speaking we have
better results for points on curves than for points on surfaces.  The
exponent for $B$ is of order $1/k$ and $1/\sqrt{k}$ in the two cases
respectively. Thus our goal is to get a result of the same exponent
$O(1/k)$ that we have for points on curves, but generalized to points
near curves.  However since ``near'' is defined here in terms of the
real metric, we will have to use a real-variable method, rather than a
$p$-adic one.

There are a number of interesting papers by Huxley \cite{hux1},
\cite{hux2}, \cite{hux3}, \cite{hux4}, relating to rational points
close to a curve.  If one
translates our problem into Huxley's language one would need to ask
about points where $f(a/q)$ is close to $b/q$, for rational $a/q$ and
$b/q$ of order 1, and an algebraic function $f$ of order 1.  It
appears that only the first of the above papers covers this particular
problem, and yields only a trivial estimate for the range of variables
which interest us.

As mentioned above, Elkies \cite{Elk} was interested in search
algorithms for rational points near curves.  It would appear that parts
of our argument can be adapted towards finding a search algorithm, but
not others.  In particular we ask: Can one find all (positive or
negative) integral solutions 
of $x_1^3+x_2^3+x_3^3=N$ in the range $\max |x_i|\le B$ in time
$O_{N,\ep}(B^{9/10+\ep})$, for any $\ep>0$? The algorithm of 
Heath-Brown \cite{HBalg} requires time $O_{N,\ep}(B^{1+\ep})$, while
that of Elkies \cite{Elk} is heuristically $O_{\ep}(B^{1+\ep})$,
uniformly for $N\le B$.

A further application of the ideas of this paper concerns the form
\[F(x_1,x_2,x_3)=x_1^{k-1}x_2-x_3^k\]
(which is singular).  In this case
an upper bound for the number of solutions of $F(\x)=N$ leads to
information on $(k-1)$-free values of the polynomial $f(x)=x^k+N$.
(An integer $n$ is said to be $l$-free if it is not divisible by the
$l$-th power of a prime.)  In fact we already know, following Hooley
\cite{hoo}, how to handle $(k-1)$-free values of a general polynomial
$f(x)$ of degree $k$, so as to give an asymptotic formula $c_fB+o(B)$
for the number of positive integers $n\le B$ with $f(n)$ being
$(k-1)$-free.  However the corresponding problem in which $n$ is
restricted to primes is much harder, and has not yet been completely solved.
The case $k\ge 7$ has been dealt with by Nair \cite[Corollary
4]{nair}, and a number of other cases are covered by Helfgott \cite{helf1}

\begin{theorem}\label{sqfr}
Let $h\in\Z-\{0\}$ be given.  Then for any $k\ge 3$ we have
\[\#\{p\le X: p^k+h\;\mbox{is}\; (k-1)\mbox{-free}\}=c_{h,k} 
{\rm Li}(X)+o(X(\log X)^{-1}),\]
where
\[c_{h,k}=\prod_{p\nmid h}(1-\frac{\nu(p)}{p^{k-2}(p-1)}),\]
with
\[\nu(r)=\#\{n\mod{r^{k-1}}: r^{k-1}|n^k+h\}.\]
\end{theorem}

In particular we can handle polynomials $x^3+h$, which in general have
Galois group $S_3$.  These are excluded from Helfgott's work
\cite{helf1}, but are covered in work to appear \cite{helf2}. We
should emphasize that our approach has little in common with
Helfgott's.  We are able to achieve substantially better error terms
than he does, but for a much restricted class of polynomials.

From now on we shall assume that the form $F$ in Theorems \ref{th9/10}
and \ref{th1/k} is fixed, and that its
degree $k$ is at least 3.  Thus all
constants implied by the $O(\ldots)$ and $\ll$ notations may depend on
$F$, and in particular on its degree $k$.

\section{The Determinant Method}\label{DM}

For the most part we will handle Theorems \ref{th9/10} and \ref{th1/k}
simultaneously.  It is only in the endgame that our arguments
diverge.  Thus we shall take $F(x_1,x_2,x_3)$ to be a fixed
non-singular integral form of degree $k$.  In this section we will
state our main lemmas and show how they lead to two of our principal
results.  Proofs of the lemmas will follow in section \ref{lemmas}.

The first move in the determinant method is to cover the region in
which points are sought via small ``patches'', which in this case will
be defined by real considerations.  We have the following lemma.
\begin{lemma}\label{patches}
There is an $M_0\in\N$, depending only on $F$, with the following properties.
Let $M$ be a positive integer. 
Suppose the square $[-1,1]^2$ is covered by $O(M^2)$ smaller
squares 
\[S=[a,a+(M_0M)^{-1}]\times[b,b+(M_0M)^{-1}].\]
Then the number 
of such squares containing a solution $(t_1,t_2)\in\R^2$ of the 
inequality $|F(t_1,t_2,1)|\le (M_0M)^{-1}$ is $O(M)$. Moreover for each
such $S$ there is an index $i=1$ or $2$ such that 
\[\left|\frac{\partial F}{\partial x_i}(a,b,1)\right|\gg 1.\] 
\end{lemma}
In the $p$-adic version of the determinant method, the analogue of
this lemma is the statement that the congruence $F(t_1,t_2,1)\equiv
0\;\;\mod{p}$ has $O_k(p)$ solutions.

We will call the squares produced by Lemma \ref{patches} ``good''.
We choose a particular good square
$S=[a,a+(M_0M)^{-1}]\times[b,b+(M_0M)^{-1}]$ and set 
\[t_1=a+u,\;\;t_2=b+v,\;\; \delta=F(t_1,t_2,1), \;\;\mbox{and}\;\;
w=\delta-F(a,b,1).\]
We then have the following result.
\begin{lemma}\label{implicit}
Suppose that $|\partial F(a,b,1)/\partial x_1|\gg 1$.  Then for 
each $s\in\N$ there are polynomials $X_s(v,w)$ and $Y_s(u,v,w)$
with the following properties:-
\begin{enumerate}
\item[(i)] The coefficients of $X_s$ and $Y_s$ are of size $O_{s}(1)$;
\item[(ii)] The total degrees of $X_s$ and $Y_s$ are bounded in
terms of $s$ and $k$;
\item[(iii)] The polynomial $X_s$ has no
constant term and $Y_s$ has no terms of total degree less than $s$;
\item[(iv)] We have $u=X_s(v,w)+uY_s(u,v,w)$.
\end{enumerate}
\end{lemma}

This lemma may be viewed as a form of the Implicit Function Theorem. It
shows that if $u$ is given implicitly by $F(a+u,b+v,1)=F(a,b,1)+w$, then
$u$ is approximately equal to $X_s(v,w)$.

We are now ready to describe the determinant method.  We begin by
choosing a positive integer $h<k$.  In fact we will later take $h=2$ in proving
Theorem~\ref{th9/10}, and $h=[(k-1)/2]$ for Theorem \ref{th1/k}. We will
consider points counted by $\cN(B;N,F,h)$ for which $B/2<\max
|x_i|=x_3\le B$. This clearly suffices, by dyadic subdivision and
symmetry.  We now write
$t_1=x_1x_3^{-1}$ and $t_2=x_2x_3^{-1}$, and proceed to apply Lemma
\ref{patches}, taking 
\begin{equation}\label{c0}
M\le (B/2)^kN^{-1}M_0^{-1}.
\end{equation}
Thus it suffices to consider separately
the relevant points in each of $O(M)$ good squares.

We choose a particular good square $S$ of side $(M_0M)^{-1}$, and 
write $(t_1^{(i)},t_2^{(i)})$, for $i\le I$ say, for the points in 
it corresponding to integer solutions of 
\begin{equation}\label{ineq}
|F(x_1^{(i)},x_2^{(i)},x_3^{(i)})|\le N.
\end{equation}
We will also
write $H=(h+1)(h+2)/2$ for convenience.  We now examine the $I\times
H$ matrix, ${\cal M}_S$ say, whose $i$-th row contains the various monomials in
$x_1^{(i)},x_2^{(i)},x_3^{(i)}$, of degree $h$.  Our aim is to
show that this matrix has rank strictly less than $H$.  We will then
be able to deduce that there are integer coefficients
$\lambda_{a,b,c}$, not all zero, such that
\[\sum_{a+b+c=h}\lambda_{a,b,c}x_1^{(i)a}x_2^{(i)b}x_3^{(i)c}=0\]
for every $i$.  This therefore produces a non-zero integral form 
$A_S(x_1,x_2,x_3)$ of degree $h$, such that $A_S(\x^{(i)})=0$ for every
$\x^{(i)}$  which corresponds to a point in $S$ and satisfies (\ref{ineq}).
We may also observe that the coefficients $\lambda_{a,b,c}$ above,
which can be constructed as certain subdeterminants of ${\cal M}_S$,
are of size $O(B^{hH})$.  In the case in which $h=2$ we deduce that
$||A_S||=O(B^{12})$, where $||A_S||$ denotes the maximum of the moduli
of the coefficients of $A_S$.

In showing that ${\cal M}_S$ has rank strictly less than $H$ we may assume
that $I\ge H$, since the assertion is trivial otherwise.  We proceed
to examine the $H\times H$ determinant, $\Delta_1$ say, arising from $H$ points
$x_1^{(i)},x_2^{(i)},x_3^{(i)}$, which without loss of generality we take to
correspond to $i=1,2,3,\ldots,H$.  By removing a factor
$x_3^{(i)h}$ from the $i$-th row, for each $i$, we find that
\begin{equation}\label{D1D2}
\Delta_1=(\prod_{i=1}^Hx_3^{(i)})^h\Delta_2\ll B^{hH}|\Delta_2|,
\end{equation}
where $\Delta_2$ is the $H\times H$ determinant whose $i$-th row
contains the monomials in $t_1^{(i)},t_2^{(i)}$ of degree at most $H$.

We shall now apply Lemma \ref{implicit}, taking
$s=H(H-1)/2$.  We have $F(a,b,1)\ll_F M^{-1}$ for any good square, and
\[\delta=F(t_1,t_2,1)=F(x_1,x_2,x_3)x_3^{-k}\ll_F NB^{-k}\ll_F
M^{-1}\]
by (\ref{c0}), whence $w\ll M^{-1}$.  Thus
$u=X_s(v,w)+O(M^{-H(H-1)/2})$, since we also
have $u,v\ll M^{-1}$.  If we now replace $w$ by $\delta-F(a,b,1)$ it 
follows that
for each monomial $t_1^et_2^f$ there is a polynomial $G_{e,f}(v,\delta)$
such that
\[t_1^et_2^f=G_{e,f}(v,\delta)+O(M^{-H(H-1)/2}).\]
The polynomial $G_{e,f}(v,\delta)$ will depend on $F$, on $S$, and on
$H$, and will have coefficients of size $O(1)$.
Now, if we write $v^{(i)},\delta^{(i)}$ for the values corresponding
to the point $t_1^{(i)},t_2^{(i)}$ and denote by $\Delta_3$ the 
$H\times H$ determinant whose
$i$-th row consists of the polynomials $G_{e,f}(v^{(i)},\delta^{(i)})$, we see
that 
\begin{equation}\label{D2D3}
\Delta_2=\Delta_3+O(M^{-H(H-1)/2}).
\end{equation}

We now employ a result relating to generalized van der
Monde determinants involving polynomials
$f_j(x_1,\ldots,x_n)\in\C(x_1,\ldots,x_n)$. We first introduce some
notation.  We let $X_1,\ldots,X_n\ge 1$ be given and we define the
size of a monomial by 
\[||x_1^{e_1}\ldots x_n^{e_n}||:=X_1^{e_1}\ldots X_n^{e_n}.\]
Writing $m$ for a typical monomial, we list the
monomials as $m_1,m_2,\ldots$ in such a way that $||m_1||\ge ||m_2||\ge
\ldots$.  Finally, as above, we shall write $||f_i||$ for the height of the
polynomial $f_i$, that is to say the maximum of the moduli of the
coefficients of $f_i$.  The result we shall use is then the following.
\begin{lemma}\label{vdM}
Let $f_1,\ldots,f_H$ be polynomials as above, having degree at most
$D$, and let
$\x^{(1)},\ldots,\x^{(H)}\in\C^n$ be vectors with
$|x_j^{(i)}|\le X_j$ for all $i$ and $j$.  Then 
\[\left|f_j(\x^{(i)})\right|_{i,j\le H}\ll_{H,D}
(\max_j||f_j||)^H\prod_{i=1}^H||m_i||.\]
\end{lemma}
 
We proceed to apply Lemma \ref{vdM} to the determinant $\Delta_3$ whose
entries are $G_{e,f}(v^{(i)},\delta^{(i)})$.  We will look at two specific 
situations.  In the first, we arrange that
the first $H$ monomials $m_i=v^e\delta^f$ are just
$1,v,v^2,\ldots,v^{H-1}$, by insisting that 
\begin{equation}\label{c1}
N(B/2)^{-k}\le M^{-(H-1)}.
\end{equation}
In this case we conclude that $\Delta_3\ll M^{-H(H-1)/2}$, and hence that 
$\Delta_2\ll M^{-H(H-1)/2}$, by (\ref{D2D3}). We deduce from
(\ref{D1D2}) that $\Delta_1\ll B^{hH}M^{-H(H-1)/2}$ and hence that
$|\Delta_1|<1$, providing that we take 
\begin{equation}\label{Mh}
M=cB^{2h/(H-1)}=cB^{4/(h+3)},
\end{equation}
with a suitably large constant $c=c(F)$. This choice is compatible
with conditions (\ref{c0}) and (\ref{c1}) providing that 
\begin{equation}\label{c2}
N\le c'B\;\;\;\mbox{and}\;\;\; h=[(k-1)/2],
\end{equation}
with a suitably small positive constant $c'=c'(F)$.
Since $\Delta_1$ has integer
entries we now conclude that $\Delta_1=0$, which establishes our
claim.  This is enough to show that the $I\times
H$ matrix ${\cal M}_S$, whose $i$-th row contains the various monomials in
$x_1^{(i)},x_2^{(i)},x_3^{(i)}$ of degree $H$, has rank strictly less than $H$.

The second situation we will examine is that in which $h=2$ (so that
$H=6$) and the first 4 monomials are $1,v,v^2,v^3$.  Thus
we assume that 
\begin{equation}\label{c3}
M^{-3}\ge N(B/2)^{-k}.
\end{equation}
The fifth and sixth monomials under our ordering can then be only
$v^4,v^5$, or $v^4,\delta$, or $\delta,v^4$.  Arguing as before we
deduce that 
\[\Delta_3\ll M^{-15}+M^{-10}NB^{-k}\ll M^{-15}+M^{-10}NB^{-3}, \]
since we are taking $k\ge 3$.  We then conclude that $\Delta_2\ll
M^{-10}NB^{-3}+M^{-15}$, and that $\Delta_1\ll
M^{-10}NB^9+M^{-15}B^{12}$.  It follows that $\Delta_1=0$ if
$M=cB^{9/10}N^{1/10}$.  This choice is compatible with (\ref{c0}) and
(\ref{c3}) when
$N\ll B^{3/13}$.  Hence, under these assumptions, we again deduce that
our $I\times H$ matrix ${\cal M}_S$ has rank strictly less than $H$.

We may now summarize our conclusions in the following two
propositions.
\begin{prop}\label{Ph}
Let $F(x_1,x_2,x_3)\in\Z[x_1,x_2,x_3]$ be a non-singular form of 
degree $k\ge 3$, and let $h$ be a positive integer with $h<k/2$.  Let
$N\in\N$ and $B\ge 1$ be given, such that $N\ll B$.  Then there are
$O(B^{4/(h+3)})$ non-zero integral forms $A_i(x_1,x_2,x_3)$ of degree
$h$, such that every integer vector with $|F(\x)|\le N$ satisfies at
least one of the equations $A_i(\x)=0$.
\end{prop}
\begin{prop}\label{P2}
Let $F(x_1,x_2,x_3)\in\Z[x_1,x_2,x_3]$ be a non-singular form of 
degree $k\ge 3$.  Let
$N\in\N$ and $B\ge 1$ be given, such that $N\ll B^{3/13}$.  Then there are
$O(B^{9/10}N^{1/10})$ non-zero integral quadratic forms
$A_i(x_1,x_2,x_3)$, such that every integer vector with 
$|F(\x)|\le N$ satisfies at
least one of the equations $A_i(\x)=0$.  Moreover the forms $A_i$ have
$||A_i||=O(B^{12})$. 
\end{prop}

In order to complete the proofs of Theorems \ref{th9/10} and
\ref{th1/k} it suffices to count points satisfying a pair of
conditions $F(\x)=N$ and $L(\x)=0$.  This will be accomplished in \S 
\ref{curves}.

\section{Lemmas \ref{patches}, \ref{implicit} and \ref{vdM}}
\label{lemmas}

For the proof the proof of Lemma \ref{patches}
it will be convenient to call a point $(t_1,t_2)\in[-1,1]^2$
satisfying $|F(t_1,t_2,1)|\le (M_0M)^{-1}$ ``good'', and similarly to call
a square containing such a point ``good''.
We begin by observing that the function
\[\max\{|\frac{\partial F}{\partial x_1}(t_1,t_2,1)|,\,
|\frac{\partial F}{\partial x_2}(t_1,t_2,1)|,\,
|\frac{\partial F}{\partial x_3}(t_1,t_2,1)|\}\]
is continuous for $(t_1,t_2)$ in the compact set $[-1,1]^2$, and is
strictly positive, since $F$ is non-singular.  Thus there is a
positive constant $\lambda$ say, depending only on $F$, such that,
for each $t_1,t_2$, at
least one of the partial derivatives
\[F_i:=\frac{\partial F}{\partial x_i}(t_1,t_2,1)\] 
has modulus $|F_i|\ge\lambda$.  Now suppose that $(t_1,t_2)$ is good,
and that $M_0\ge 3k/\lambda$.  By Euler's identity we have
\[|F_1t_1+F_2t_2+F_3|=|kF(t_1,t_2,1)|\le k/(M_0M)\le \lambda/3.\]
Thus if $|F_3|\ge \lambda$ we must have $|F_it_i|\ge\lambda/3$ for
either $i=1$ or 2, and hence $|F_i|\ge\lambda/3$.  It follows that we
will necessarily have $|F_i|\ge\lambda/3$ for either $i=1$ or $i=2$,
or both.

Since the partial
derivatives are continuous, there is an integer $M_0$, depending only
on $F$, such that $F_1$ and $F_2$ vary by at most $\lambda/6$ over any
square $S_0\subseteq [-1,1]^2$ of side $M_0^{-1}$.  Each good point
therefore
lies in a square $S_0$ of side $M_0^{-1}$ such that for some choice of
$i=1$ or 2, and some choice of $\pm$ sign, we have 
\[\pm F_i(x_1,x_2,1)\ge \lambda/6,\;\;\;\mbox{for all}\;\;\;
(x_1,x_2)\in S_0.\]

Let us consider a square $S_0$ for which the index $i=1$ and the $+$
sign are admissible.  We proceed to cover $S_0$ with $M^2$ 
squares, of the type
\[S_{u,v}:=\left[a+\frac{u-1}{MM_0}\,,\,a+\frac{u}{MM_0}\right]\times
\left[b+\frac{v-1}{MM_0}\,,\,b+\frac{v}{MM_0}\right],\;\;\;
(1\le u,v\le M).\]
Suppose that there are two good squares $S_{u,v}$ and $S_{u',v}$.
Then, by the Mean-Value Theorem, there is a point $(\xi_1,\xi_2)$ on
the line between $(t_1,t_2)$ and $(t_1',t_2')$, such that
\[|(t_1-t_1')F_1(\xi_1,\xi_2,1)+(t_2-t_2')F_2(\xi_1,\xi_2,1)|=
|F(t_1,t_2,1)-F(t_1',t_2,1')|\le 2M^{-1}.\]
Taking 
\[M_0\ge\sup \{|F_2(x_1,x_2,1)|:\,(x_1,x_2)\in[-1,1]^2\},\]
as we may, we deduce that $|(t_1-t_1')F_1(\xi_1,\xi_2,1)|\le 3M^{-1}$,
and hence that $|t_1-t_1'|\le 18(M\lambda)^{-1}$.  It follows that
$S_0$ contains $O(1)$ good squares $S_{u,v}$, for each fixed $v$.

Finally we deduce that $S_0$ contains $O(M)$ good squares $S_{u,v}$,
and hence that we can cover $[-1,1]$ with $O(M)$ good squares $S_{u,v}$, of
side $(MM_0)^{-1}$, as required.
\bigskip

We turn now to the proof of Lemma \ref{implicit}.
We begin by observing that $w=uF_1+vF_2+f(u,v)$ for some
polynomial $f$ composed of monomials of degree at least 2, where $F_1$
and $F_2$ are the usual partial derivatives at $(a,b,1)$.  Moreover
it is clear that the coefficients of $f$ are $O(1)$. Since $F_1\gg 1$ 
we may rewrite the equation in the form
\begin{eqnarray*}
u&=&\{w-vF_2-f(0,v)\}F_1^{-1}+u\{(f(0,v)-f(u,v))/u\}F_1^{-1}\\
&=&X_1(v,w)+uY_1(u,v,w)
\end{eqnarray*}
say, where $X_1$ and $Y_1$ satisfy (i), (ii) and
(iii).  We will
deduce by an inductive iteration that for any $s\ge 1$ we may write 
$u=X_s(v,w)+uY_s(u,v,w)$, where $X_s$ and $Y_s$ satisfy the conditions
of the lemma.  Specifically we have
\begin{eqnarray*}
u&=&X_s(v,w)+uY_s(u,v,w)\\
&=&X_s(v,w)+\{X_s(v,w)+uY_s(u,v,w)\}Y_s(u,v,w)\\
&=&X_s(v,w)+X_s(v,w)Y_s\{X_s(v,w)+uY_s(u,v,w),v,w\}
+uY_s(u,v,w)^2\\
&=&X_{s+1}(v,w)+uY_{s+1}(u,v,w),
\end{eqnarray*}
where
\[X_{s+1}(v,w)=X_s(v,w)\{1+Y_s\left(X_s(v,w),v,w\right)\},\]
and
\[Y_{s+1}(u,v,w)=X_s(v,w)\frac{Y^{(1)}-Y^{(2)}}{u}+Y_s(u,v,w)^2,\]
where
\[Y^{(1)}=Y_s\{X_s(v,w)+uY_s(u,v,w),v,w\}\;\;\;\mbox{and}\;\;\;
Y^{(2)}=Y_s\{X_s(v,w),v,w\}.\]
One can easily verify that $X_{s+1}$ and $Y_{s+1}$ have the required
properties, and the induction is complete.
\bigskip

We end this section by considering Lemma \ref{vdM}.  We shall assume that
\[\max_j||f_j||\ll_H 1,\]
as we clearly may.  Recall that we have given the monomials in
$x_1,\ldots,x_H$ an ordering such that $||m_1||\ge ||m_2||\ge \ldots.$
We shall call $m_i$ the ``leading monomial'' in a
polynomial $f$, if it is the monomial with non-zero coefficient for
which $i$ is least.  We shall then say that $i$ is the
``index'' of $f$, and write ${\rm ind}(f)=i$.

We shall perform a sequence of at most $H^2$ 
elementary column operations on the determinant
\[\Delta:={\rm det}\left(f_j(\x^{(i)})\right),\]
which will replace the set of polynomials $f_1,\ldots,f_H$ by a new set 
$g_1,\ldots,g_H$, so that
\[\pm\Delta=\Delta':={\rm det}\left(g_j(\x^{(i)})\right),\]
The new polynomials will have degree at most $D$, and will also satisfy
\[\max_j||g_j||\ll_H 1.\]
Moreover they will have the property that
\[{\rm ind}(g_1)<{\rm ind}(g_2)<\ldots<{\rm ind}(g_H),\]
whence ${\rm ind}(g_j)\ge j$. Now for any polynomial
$g(\x)\in\C[x_1,\ldots,x_n]$ of index $r$ we have
\[|g(\x)|\ll_D ||g||.||m_r||.\]
It therefore follows that the $ij$ entry in $\Delta'$ is $O_{H,D}(||m_{j}||)$,
whence
\[\Delta'\ll_{H,D}\prod_{j=1}^H||m_j||.\]
The lemma then follows, once we have shown how the polynomials $g_j$
are obtained.

We shall show by induction on $r$ that, after at most $rH$ column operations, 
we can ensure that we have polynomials with 
${\rm ind}(g_1)<{\rm ind}(g_2)<\ldots<{\rm ind}(g_r)$ and 
${\rm ind}(g_j)>{\rm ind}(g_n)$ for $j>r$.  When $r=H$ this gives the
required result.  The base case of the induction, in which $r=0$, is
trivial, so we shall assume that the above statement holds for $r=s-1$
say, and prove that it also holds for $r=s$.   Suppose that
$i$ is the smallest index occurring among the polynomials $g_s,\ldots,
g_H$.  Of all such polynomials with index $i$ we choose one, $g_j$
say, for which the coefficient of $m_i$ is largest in
modulus, and swap the $j$-th and $s$-th columns of the
determinant. After this re-ordering of the polynomials we will have
\[{\rm ind}(g_1)<{\rm ind}(g_2)<\ldots<{\rm ind}(g_s)=i\]
and 
${\rm ind}(g_j)\ge i={\rm ind}(g_s)$ for $j>s$.
For each $j=s,\ldots,H$ we now write $c_j$ for the coefficient of
$m_i$ in $g_j$, so that $|c_j|\le |c_s|$ for $j\ge s$.
We proceed to perform further column operations on the determinant, 
subtracting $c_jc_s^{-1}$ times column $s$ from column $j$, for
$j>s$.  This produces a determinant with new polynomials $g_j$ for $j>s$,
satisfying ${\rm ind}(g_j)>i={\rm ind}(g_s)$.  This establishes the
induction step, since we have used at most $H$ column operations.

\section{Counting Points on Curves}\label{curves}

To complete the proofs of Theorems \ref{th9/10} and \ref{th1/k} we
need to estimate the number of points on the affine curves
\[F(\x)=N,\;\;\;A_i(\x)=0,\]
where $A_i$ has degree $h=2$ or $[(k-1)/2]$.
We shall use techniques from the author's work \cite[\S 9]{annals}. 

We first consider Theorem \ref{th1/k}.  The equations
\[F(x_1,x_2,x_3)-Nx_4^k=A_i(x_1,x_2,x_3)=0\]
define a projective curve in $\Proj^3$, of degree $hk$.  This curve
may or may not be
irreducible.  According to \cite[Theorem 5]{annals}, any component of
degree $d$ will contribute $O_{\ep}(B^{2/d+\ep})$ to 
$\cN(B;N,F,h)$.  In particular, if $d\ge k-1$, we get a total
contribution $O_{\ep}(B^{\phi+\ep})$ on considering the various
possible forms $A_i$, where 
\[\phi=\frac{2}{k-1}+\frac{4}{h+3}\le\frac{2}{k-1}+\frac{4}{(k-2)/2+3}
=\frac{2}{k-1}+\frac{8}{k+4}<\frac{10}{k}.\]
Thus components of degree at least $k-1$ make an acceptable
contribution.  According to a theorem of Colliot-Th\'{e}l\`{e}ne, see
\cite[Appendix]{annals}, the surface 
\[F(x_1,x_2,x_3)-Nx_4^k=0\]
contains $O_k(1)$ curves of degree at most $k-2$, since $F$ is
non-singular.  The projection of such a curve onto the plane $x_4=0$
produces a projective plane curve $J(x_1,x_2,x_3)=0$.  Here $J$ will
of course be a factor of one of the forms $A_i$.  Thus we have to
consider $O_k(1)$ intersections
\begin{equation}\label{int}
F(x_1,x_2,x_3)=N,\;\;\;J_i(x_1,x_2,x_3)=0,
\end{equation}
where the set of available forms $J_i$ is independent of the choice of
$N$, as in \cite[\S 9]{annals}.  When $J_i=0$ defines a plane curve of
genus at least 2, it will have finitely many projective points, and in
particular the number of such points is independent of $N$, but not
necessarily independent of $F$.  Each
projective point produces at most two solutions of $F(\x)=N$, so that
the overall contribution to $\cN(B;N,F,h)$ is $O_F(1)$.  When $J_i=0$
has genus 1 there will be $O_{\ep}(B^{\ep})$ projective points, with a
constant depending on $J_i$ as well as $\ep$.  Thus this case
contributes $O_{F,\ep}(B^{\ep})=O(B^{10/k})$ to $\cN(B;N,F,h)$.

When $J_i=0$ has genus zero, it can be parameterized by coprime forms 
$f_1(u,v),f_2(u,v),f_3(u,v)\in\Z[u,v]$.  We then find, as in
\cite[pp 592 \& 593]{annals}, that the solutions of (\ref{int})
are given, with $O_k(1)$ exceptions, by
\[(x_1,x_2,x_3)=(\lambda\nu^{-1}f_1(u,v),\lambda\nu^{-1}f_2(u,v),
\lambda\nu^{-1}f_3(u,v)).\]
Here $\lambda$ and $\nu$ are coprime integers, as are $u$ and
$v$. Moreover the forms $f_1,f_2,f_3$ may be assumed to be coprime.
Since $\nu$ divides each of $f_1(u,v), f_2(u,v)$ and $f_3(u,v)$ it
must also divide their resolvent, which must be non-zero, because the
forms are coprime.  Thus $\nu$ can take at most
$O_k(1)$ values in total.  The pairs $u,v$ for which $\nu$ divides
each of the $f_j(u,v)$ lie on one of $O_k(1)$ lattices.  Hence, by
making appropriate linear substitutions, we conclude that the
solutions to (\ref{int}) may be given with $O_k(1)$ exceptions, by one
of $O_k(1)$ parameterizations
\[(x_1,x_2,x_3)=(\lambda g_1(s,t),\lambda g_2(s,t),\lambda g_3(s,t)).\]
Here $s,t$ will take integer values.  Moreover $\lambda$ can take only
$O_{\ep}(N^{\ep})$ values, since $\lambda^k|N$.  
Finally we may absorb $\lambda$
into the forms $g_j$, and restrict attention to solutions with
$x_j=g_j(s,t)$. 

Now, again as in \cite[p 593]{annals}, we see that there are 
$O_{\ep}((BN)^{\ep})$
solutions $s,t$ unless $g(s,t):=F(g_1(s,t),g_2(s,t),g_3(s,t))$ is a
power of a linear form.  By an invertible linear change of variable we
may then suppose that $g(s,t)=ct^{dk}$ say, where $d$ is the degree of
the forms $g_i$ and the coefficient $c$ is a non-zero integer.  
It then follows that $N$ must take the form $N=cN_0^k$, and that
$t^d=\pm N_0$.  We now recall that the available forms
$J_i$ are determined by $F$ only, and hence so are the forms $g_j$ (up
to multiplication by $\lambda$).
Thus the solutions we must count fall into $O(1)$ parametric families,
as claimed in theorem \ref{th1/k}.  
The forms $g_j$ are determined by $F$, but are
independent of
$N$, and clearly a parameterization of degree $d$ can produce only
$O(B^{1/d})$ solutions.  Thus parameterizations of degree greater than
$[k/10]$ contribute an acceptable amount, and Theorem \ref{th1/k}
follows.

We turn now to Theorem \ref{th9/10}.  We have to consider
$O(B^{9/10}N^{1/10})$ curves
\begin{equation}\label{curve}
F(\x)=N,\;\;\;A_i(\x)=0,
\end{equation}
where $A_i$ has degree $h=2$, and $||A_i||=O(B^{12})$.  The form 
$A_i$ might be reducible, in which case we replace it by its
factors. If $A_i$ is reducible only over a quadratic extension of $\Q$
then $\x$ must satisfy a pair of linear equations. 
Thus we may assume that $A_i$ is irreducible of degree 1 or 2.  

We shall first consider the situation where $A_i$ is of degree 2.
In this case we can parameterize the solutions of $A_i(\x)$ via quadratic
forms 
\[f_1(u,v),f_2(u,v),f_3(u,v)\in\Z[u,v]\]
such that any
solution $\x$ must be proportional to $(f_1(u,v),f_2(u,v),f_3(u,v))$
for some coprime integers $u,v$.   The forms $f_j$ will be pairwise
coprime as polynomials and will have heights
$||f_j||=O(B^c)$ for some absolute constant $c$.
We then see that
\begin{equation}\label{xf}
(x_1,x_2,x_3)=\lambda\nu^{-1}(f_1(u,v),f_2(u,v),f_3(u,v)),
\end{equation}
where $\nu$ divides each of $f_1(u,v),f_2(u,v)$ and $f_3(u,v)$.  For a
fixed index $i$ the
values of $\lambda$ and $\nu$ may be different for different vectors $\x$.
However, as we took $u$ and $v$ to be coprime we see that
$\nu$ divides the resolvent ${\rm Res}(f_1,f_2)=R$, say.  Since
$R$ must be a non-zero integer of size $O(B^{4c})$, we deduce that
there are $O_{\ep}(B^{\ep})$ possible values of $\nu$, for any
$\ep>0$, for each index $i$.
We also have 
\[\lambda^kF(f_1(u,v),f_2(u,v),f_3(u,v))=\nu^k N,\]
whence $\lambda^k|N$.  Thus $\lambda$ also takes $O_{\ep}(B^{\ep})$ 
possible values, for each index $i$.

It remains to consider how many solutions an equation of the shape
\[F(f_1(u,v),f_2(u,v),f_3(u,v))=N_0\]
can have.  We shall write $G(u,v)$ for the form on the left, so that
$G$ has degree $2k$. If $G$ is a product of two coprime integral factors $G_1$
and $G_2$ say, then we must have $G_1(u,v)=N_1$ and $G_2(u,v)=N_2$ for
some pair of integers $N_1N_2=N_0$.  Since $G_1$ and $G_2$ are coprime
these equations determine $O(1)$ values of $u,v$, by elimination.  We
therefore obtain a total of $O_{\ep}(B^{\ep})$ solutions of this type, 
for a given index $i$.  The alternative is that $G=cG_1^r$ for some
integral form $G_1$, irreducible over $\Q$.  In this case we obtain a 
Thue equation of the form $G_1(u,v)=N_1$, with $u,v$ coprime.  It has
been shown by Bombieri and Schmidt \cite{thue} that if $G_1$ has
degree $2k/r\ge 3$ then the number of solutions is
$O((2k/r)^{\omega(N_1)})=O_{\ep}(B^{\ep})$, where $\omega(n)$ denotes
the number of distinct prime factors of $n$.  This case is therefore
satisfactory too.

We next examine the situation in which $G_1$ has degree 2.  Since the
coefficients of $G_1$ are bounded by powers of $B$, the theory of the
Pell equation shows that the number of solutions of
$G_1(u,v)=N_1$ will be $O_{\ep}(B^{\ep})$, providing that $u$ and $v$
are also bounded by powers of $B$.  Now it is clear from (\ref{xf})
and from what we have said about $\mu$ and $\nu$ that each of 
$f_i(u,v)$ is bounded by a power of $B$. Let
\[f_1(u,v)=a_1u^2+b_1uv+c_1v^2=\tau_1,\;\;\;\mbox{and}
\;\;\; f_2(u,v)=a_2u^2+b_2uv+c_2v^2=\tau_2.\]
Since $f_1$ and $f_2$ are coprime as forms we have
$\tau_1=\tau_2=0$ only when $u=v=0$. Otherwise we note that $u-vj$ divides
$f_2(j,1)\tau_1-f_1(j,1)\tau_2$. We may choose $j\ll 1$ so that this
last quantity is non-zero, and deduce that $u-jv$ is bounded by a
power of $B$.  Since each $f_i(u,v)$ is also bounded by a power of $B$
we may deduce that both $u$ and $v$ are bounded by powers of $B$, as
required. 

There remains the case in which $G_1$ is linear, so that we have an
identity
\[F(f_1(u,v),f_2(u,v),f_3(u,v))=cG_1(u,v)^{2k}.\]
By making a linear change of variables, invertible over $\Z$, we can
assume that $G_1(u,v)=v$, so that
\[F(f_1(u,v),f_2(u,v),f_3(u,v))=cv^{2k}.\]
It follows that the solutions, in $\Proj^2(\Qb)$, of $A_i(\x)=0$, can be
parameterized by quadratic forms $f_i$ satisfying 
\[F(f_1(u,v),f_2(u,v),f_3(u,v))=v^{2k}.\]
We shall say, under such circumstances, that $A_i$ is ``special'' for
$F$.

We can also argue as above when $A_i$ is linear.  We find that
(\ref{curve}) has $O(B^{\ep})$ relevant solutions, unless
$A_i(\x)=0$ can be
parameterized by linear forms $f_i$ satisfying 
\begin{equation}\label{nt}
F(f_1(u,v),f_2(u,v),f_3(u,v))=v^{k}.
\end{equation}
Again we shall say in this case that $A_i$ is ``special'' for
$F$.

We now call on the following lemma.
\begin{lemma}\label{qb}
Let $F(X_1,X_2,X_3)\in\Qb[X_1,X_2,X_3]$ be a non-singular form of 
degree $k\ge 3$, then the number of linear or quadratic forms $A_i$,
up to scalar multiplication, which are special for $F$, is bounded in
terms of $k$ alone.
\end{lemma}

We will prove this later, but first we show how it suffices to
complete the proof of Theorem \ref{th9/10}.  The number of forms $A_i$
to be considered is $O(B^{9/10}N^{1/10})$ and each one contributes
$O(B^{\ep})$, unless it is special.  By Lemma \ref{qb} there are
$O(1)$ special forms $A_i$.  Moreover it is clear from the definition,
that the list of special forms does not depend on $N$.  The special
forms $A_i$ of degree 1 are precisely those corresponding to linear
parameterizations as described in the introduction.  When
$F(\x)=x_1^k\pm x_2^k\pm x_3^k$, an easy argument as in \cite[p
  ???]{annals}, based on Fermat's Last Theorem, shows that
any linear parameterization must correspond only to ``special solutions''.

Hence it remains to examine the case in which $A_i$ is 
quadratic.  Here the form $A_i$ determines $\lambda$ and $\nu$, which are
therefore of size $O(1)$.  Then $N,\lambda$ and $\nu$ determine $v$ up
to sign (having arranged as above that $G_1(u,v)=v$).  For any
quadratic polynomial $f(u)\in\Z[u]$ there can be at most $O(B^{1/2})$
integers $u$ with $f(u)\ll B$, and thus each $v$ corresponds to at
most $O(B^{1/2})$ choices for $u$.  We therefore deduce that each
special quadratic $A_i$ contributes $O(B^{1/2})$ in Theorem
\ref{th9/10}, which is satisfactory.

We have now established Theorem \ref{th9/10} completely, except for
the proof of Lemma \ref{qb}.  For the latter we begin by observing
that the identity (\ref{nt}) corresponds to a line (if $A_i$ is
linear), or a conic (if $A_i$ is quadratic), lying in
the non-singular projective surface
\[F(X_1,X_2,X_3)-X_4^k=0.\]
A theorem of Colliot-Th\'{e}l\`{e}ne \cite[Appendix]{annals} shows
that the surface contains at most $O_k(1)$ curves of degree $\le
k-2$. Moreover such a curve, parameterized by linear or quadratic
forms, will determine $A_i$ by projection onto $X_4=0$.  This
establishes the lemma, except for the case of quadratic
parameterizations when $k=3$.

We now handle this final case.  The condition that $A=A_i$ is
special implies that the intersection
\[F(X_1,X_2,X_3)=A(X_1,X_2,X_3)=0\]
has only a single point, which we take to be $P=(1,0,0)$.  Hence, after a
change of variables which fixes $P$ we may suppose that
$A=X_1X_2-X_3^2$. It then follows that $F(u^2,uv,v^2)=0$ only for
$v=0$, whence $F(u^2,uv,v^2)=cu^6$ for some constant $c$.  We deduce
that $F$ takes the shape 
\[F(X_1,X_2,X_3)=A(X_1,X_2,X_3)L(X_1,X_2,X_3)+cX_2^3\]
for some linear form $L$.  Moreover the line $X_2=0$ is tangent to the
conic $A=0$ at $P$.  

In general we may conclude that if $A$ is
special then $F=AL+L'^3$ for appropriate linear forms $L$ and $L'$,
such that $L=0$ is tangent to $A=0$ at the point $Q$ where $L=L'=0$.
(Note that $L$ and $L'$ cannot be
proportional, since $F$ is non-singular.)
We now observe that $F=0$ meets the line $L=0$ only at $Q$.
Thus $Q$ must be one of the
flex points of $F=0$, and $L=0$ must be the corresponding tangent
line.  

It will therefore suffice to show that there are a finite
number of possible quadratic forms $A$ corresponding to a given pair
$Q,L$. In order to do this we may change variables so that $L=X_3$ and
$Q=(1,0,0)$.  Suppose now that 
\[F(X_1,X_2,X_3)=A_0(X_1,X_2,X_3)X_3+L_0(X_1,X_2,X_3)^3\]
for a particular special form $A_0$.  We may change variables so that
$L_0=X_2$. We then consider other possible special quadratics $A$
corresponding to the same point $Q$.  For these we have
\[F(X_1,X_2,X_3)=A(X_1,X_2,X_3)X_3+L(X_1,X_2,X_3)^3.\]
Clearly $X_3|L^3-X_2^3$, and we may choose $L$ so that $X_3$ divides
$L-X_2$.  On setting $L=X_2+tX_3$ we conclude that 
\[A_0(X_1,X_2,X_3)X_3+X_2^3=A(X_1,X_2,X_3)X_3+(X_2+tX_3)^3\]
identically.  We proceed to substitute $X_2=-tX_3$, whence
\[A_0(X_1,-tX_3,X_3)-t^3X_3^2=A(X_1,-tX_3,X_3)\]
identically in $X_1$ and $X_3$.  However, since the line $L=0$ is
tangent to the conic $A=0$ we deduce that the binary quadratic form
$A(X_1,-tX_3,X_3)$ must be a square.  It follows that if the value $t$
corresponds to a special form $A$ then $A_0(X_1,-tX_3,X_3)-t^3X_3^2$
must be a square. The determinant of this quadratic form is a
non-zero polynomial in $t$ of degree at most 3, so we may conclude
finally that there are at most 3 special quadratics corresponding to a
given flex point Q.  This completes the proof of Lemma \ref{qb}.

\section{Theorem \ref{sqfr}}

In this section we shall consider Theorem \ref{sqfr}.  We shall be
relatively brief, since the main ideas have all been expounded above.
We recall that $h$ and $k$ are considered fixed. 

We have
\begin{equation}\label{est1}
\#\{p\le X: p^k+h\;\mbox{is}\;(k-1)\mbox{-free}\}=\sum_{d=1}^{\infty}\mu(d)
\#\{p\le X: d^{k-1}|p^k+h\}.
\end{equation}
Only values of $d$ with $d^{k-1}\ll X^k$ can contribute.  Moreover
pairs $d,h$ with a non-trivial common factor produce only 
$O(1)$.  The function $\nu(r)$
in Theorem \ref{sqfr} satisfies $\nu(r)\ll_{\ep}r^{\ep}$ for any
$\ep>0$, whence the Siegel Walfisz Theorem yields
\[\#\{p\le X:
d^{k-1}|p^k+h\}=\frac{\nu(d)}{\phi(d^{k-1})}{\rm Li}(X)+O(X(\log X)^{-5})\]
for $d\le(\log X)^3$, say.  Thus
\begin{eqnarray*}
\lefteqn{\sum_{d\le (\log X)^3}\mu(d)\#\{p\le X: d^{k-1}|p^k+h\}}\hspace{2cm}\\
&=&\sum_{d\le (\log X)^3}\mu(d)\frac{\nu(d)}{\phi(d^{k-1})}{\rm Li}(X)
+O(X(\log X)^{-2})\\
&=& c_{h,k} {\rm Li}(X)+O(X(\log X)^{-2}).
\end{eqnarray*}

Using the bound $\nu(r)\ll_{\ep} r^{\ep}$ once more we have
\[\#\{p\le X: d^{k-1}|p^k+h\}\ll_{\ep} d^{\ep}(Xd^{1-k}+1).\]
Thus the range $(\log X)^3\le d\le X^{1-\ep}$ contributes $o(X(\log
X)^{-1})$ to (\ref{est1}), providing that $\ep<1$.
We handle the remaining range via the following lemma.
\begin{lemma}\label{sfl}
Let $k\ge 3$, and suppose that $A$ and $B$ are real parameters with 
$B\ge A^{1-1/(4k+3)}$.
Then the number of integer solutions $x,y,z$ of the
equation $x^k+h=y^{k-1}z$ with $A<x\le 2A$ and $B<y\le 2B$ will be
$O_{\ep}(A^{19/20}B^{\ep})$, for any fixed $\ep>0$.
\end{lemma}
Clearly this suffices for the proof of Theorem \ref{sqfr}.  Note that
this estimate is strongly related to that given by Theorem
\ref{th9/10}, but that in the present case the form
$F(x,y,z)=y^{k-1}z-x^k$ is singular.

Following the strategy for the proof
of Lemma \ref{patches} we break the available
range $[A/(2B),2A/B]$ for $x/y$ into $M$ equal subintervals of the
form $[a,a+3A/(2BM)]$, and write $x/y=a+v$ with $0\le v\le 3A/(2BM)$.
Then 
\[z/y=(x/y)^k+hy^{-k}=(a+v)^k+\delta,\]
say, with $\delta\ll B^{-k}\ll B^{-3}$. This identity plays the role
of Lemma \ref{implicit}.  We proceed to form the matrix of quadratic
monomials corresponding to solutions $x,y,z$ with $x/y$ in a given
interval $[a,a+3A/(2BM)]$.  The argument then continues as in \S
\ref{DM}, with the simplification that there is no error term
$O(M^{-H(H-1)/2})$ in (\ref{D2D3}) to deal with.  The condition 
corresponding to 
(\ref{c3}) is that $(A/BM)^3\gg B^{-3}$, or merely that $M\ll A$.
When we apply Lemma
\ref{vdM} we have to allow for the fact that our polynomials depend on
$a$, which may not be of order 1.  Thus, assuming that $M\ll A$, we
will get
\begin{eqnarray*}
\Delta_3&\ll & 
\max(1\,,\,|a|)^{2k}\left\{(A/BM)^{15}+(A/BM)^{10}B^{-3}\right\}\\
&\ll & (1+A/B)^{2k}\left\{(A/BM)^{15}+(A/BM)^{10}B^{-3}\right\}.
\end{eqnarray*}
We therefore conclude that
\[\Delta_1\ll B^{12}(1+A/B)^{2k}\left\{(A/BM)^{15}+(A/BM)^{10}B^{-3}\right\}.\]
This shows that $|\Delta_1|<1$ if 
\[M\gg \frac{A}{B^{1/10}}\left\{1+(\frac{A}{B})^{2k/10}\right\}.\]
In particular it suffices to have $M\ge A^{19/20}$ if $B\ge A^{1-1/(4k+3)}$
and $A$ is sufficiently large.

It follows that the solutions to $x^k+h=y^{k-1}z$ with $A<x\le 2A$ 
and $B<y\le 2B$ will lie on $O(A^{19/20})$ lines or conics
$A_i(x,y,z)=0$.  The argument for \S \ref{curves} goes through just as
before, until we reach Lemma \ref{qb}, at which point crucial use was
made of the fact that $F$ was non-singular. Thus it remains to
consider whether there might exist non-zero linear or quadratic polynomials
$f_1(u), f_2(u)$ and $f_3(u)$, not all constant, such that
\[f_1(u)^k-f_2(u)^{k-1}f_3(u)=1.\]
Since $k\ge 3$ it is clear that neither $f_1$ nor $f_2$ can 
be constant.  If $l(u)$ is a linear factor of $f_2(u)$ then
\[l(u)^{k-1}|f_1(u)^k-1=\prod_{i=0}^{k-1} (f_1(u)-\omega^i),\]
where $\omega$ is a primitive $k$-th root of unity.  Since the factors
on the right differ by a non-zero constant they must be coprime,
whence $l(u)^{k-1}|f_1(u)-\omega^i$ for some $i$.  This is only
possible if $k=3$ and $f_1$ is quadratic, whence $f_2$ and $f_3$ must
also be quadratic.  We then see that $f_2(u)=l_1(u)l_2(u)$ with
$l_1(u)^{k-1}|f_1(u)-\omega^i$ and $l_2(u)^{k-1}|f_1(u)-\omega^j$ for
two different exponent $i,j$.  Thus both $f_1(u)-\omega^i$ and
$f_1(u)-\omega^j$ must be squares, say $f_1(u)-\omega^i=l(u)^2$ and
$f_1(u)-\omega^j=l'(u)^2$.  However this is impossible since
$l(u)^2-l'(u)^2=(l(u)-l'(u))(l(u)+l'(u))$ cannot be a
non-zero constant.  We therefore conclude that "special" forms $A_i$
(in the sense of Lemma \ref{qb}) do not exist in the current setting.
We deduce that each line or conic $A_i(x,y,z)=0$ produces
$O_{\ep}(B^{\ep})$ solutions, and Lemma \ref{sfl} follows.

\bigskip
\bigskip

Mathematical Institute,

24--29, St. Giles',

Oxford

OX1 3LB

UK
\bigskip

{\tt rhb@maths.ox.ac.uk}

\end{document}